\documentclass[]{article}          \usepackage{authblk}
\usepackage[margin=2cm]{geometry}

\RequirePackage{calrsfs}
\RequirePackage{amsmath}
\RequirePackage{amsthm}

\RequirePackage{amssymb}
\RequirePackage{amsfonts}

\RequirePackage[most]{tcolorbox}

\RequirePackage{bbm}

\RequirePackage{graphicx}

\RequirePackage{booktabs}

\RequirePackage{url}

\RequirePackage{hyperref}

\RequirePackage{xcolor}
\RequirePackage{color}

\RequirePackage{wasysym}
\RequirePackage{listings}
\RequirePackage{caption}
\RequirePackage{subcaption}

\RequirePackage{stmaryrd}

\RequirePackage{lineno}

\RequirePackage{xparse}

\RequirePackage{float}

\RequirePackage{multirow}

\RequirePackage{cellspace}

\RequirePackage{shellesc}

\RequirePackage[ruled,vlined,linesnumbered]{algorithm2e}

\SetCommentSty{xCommentSty}

\RequirePackage{tikz}
\usetikzlibrary{calc, automata, positioning, arrows, external, matrix, shapes.misc}
\tikzset{cross/.style={cross out,thick,draw=black,minimum size=2*(#1-\pgflinewidth), inner sep=0pt, outer sep=0pt},
cross/.default={2.5pt}}
  
\RequirePackage{pgfplots}
\pgfplotsset{compat=1.16}
\usepgfplotslibrary{fillbetween}
\usepgfplotslibrary{groupplots}
\usepgfplotslibrary{external}

\DeclareMathAlphabet{\pazocal}{OMS}{zplm}{m}{n}
\DeclareMathAlphabet\bpazocal{OMS}{cmsy}{b}{n}

\providecommand{\Ker}{\ensuremath{\text{Ker}}}

\providecommand{\bR}{\mathbb{R}}
\providecommand{\bs}{\boldsymbol}

\providecommand{\T}{\mathsf{T}}

\providecommand{\half}{\frac{1}{2}}

\providecommand{\vertiii}[1]{{\left\vert\kern-0.15ex\left\vert\kern-0.15ex\left\vert #1
    \right\vert\kern-0.15ex\right\vert\kern-0.15ex\right\vert}}

\NewDocumentCommand{\curlii}{sO{}m}
{
	\IfBooleanTF{#1}
    {\dgalext{#3}}
    {\dgalx[#2]{#3}}
}

\NewDocumentCommand{\dgalext}{m}{  \sbox0{    \mathsurround=0pt     $\left\{\vphantom{#1}\right.\kern-\nulldelimiterspace$  }  \sbox2{\{}  \ifdim\ht0=\ht2
    \{\kern-.625\wd2 \{#1\}\kern-.625\wd2 \}  \else
    \left\{\kern-.7\wd0\left\{#1\right\}\kern-.7\wd0\right\}  \fi
}

\NewDocumentCommand{\dgalx}{om}{  \sbox0{\mathsurround=0pt$#1\{$}  \sbox2{\{}  \ifdim\ht0=\ht2
    \{\kern-.625\wd2 \{#2\}\kern-.625\wd2 \}  \else
    \mathopen{#1\{\kern-.7\wd0 #1\{}
    #2
    \mathclose{#1\}\kern-.7\wd0 #1\}}
  \fi
}

\tikzset{
  partial ellipse/.style args={#1:#2:#3}{
    insert path={+ (#1:#3) arc (#1:#2:#3)}
  }
}

\providecommand{\pE}{\pazocal{E}}

\providecommand{\pG}{\pazocal{G}}

\providecommand{\pO}{\pazocal{O}}

\providecommand{\pQ}{\pazocal{Q}}

\providecommand{\pT}{\pazocal{T}}

\providecommand{\bpM}{\bpazocal{M}}

\providecommand{\bpV}{\bpazocal{V}}

\providecommand{\matA}{\bs{A}}
\providecommand{\matB}{\bs{B}}

\providecommand{\matD}{\bs{D}}

\providecommand{\matH}{\bs{H}}
\providecommand{\matI}{\bs{I}}

\providecommand{\matP}{\bs{P}}

\providecommand{\matS}{\bs{S}}

\providecommand{\matX}{\bs{X}}

\providecommand{\vecb}{\bs{b}}

\providecommand{\vece}{\bs{e}}
\providecommand{\vecf}{\bs{f}}

\providecommand{\vecn}{\bs{n}}

\providecommand{\vect}{\bs{t}}
\providecommand{\vecu}{\bs{u}}
\providecommand{\vecv}{\bs{v}}

\providecommand{\vecx}{\bs{x}}
\providecommand{\vecy}{\bs{y}}

\providecommand{\to}{\widetilde{o}}

 \graphicspath{{./pngfigs/}}
\providecommand{\keywords}[1]
{
  \small	
  \textbf{\textit{\ \ \ \ Keywords:}} #1
}

\makeatletter
\if@twocolumn

\else

\fi
\makeatother

\begin{document}

\title{Multigrid and saddle-point preconditioners for unfitted finite element modelling of inclusions}

\author[]{Hardik Kothari\thanks{hardik.kothari@usi.ch}}
\author[]{Rolf Krause\thanks{rolf.krause@usi.ch}}
\affil[]{Institute of Computational Science, Universit\`{a} della Svizzera italiana, Lugano, Switzerland}

\maketitle

\begin{abstract}
  In this work, we consider the modeling of inclusions in the material using an unfitted finite element method.
In the unfitted methods, structured background meshes are used and only the underlying finite element space is modified to incorporate the discontinuities, such as inclusions.
Hence, the unfitted methods provide a more flexible framework for modeling the materials with multiple inclusions.
We employ the method of Lagrange multipliers for enforcing the interface conditions between the inclusions and matrix, this  gives rise to the linear system of equations of saddle point type.
We utilize the Uzawa method for solving the saddle point system and propose preconditioning strategies for primal and dual systems.

For the dual systems, we review and compare the preconditioning strategies that are developed for FETI and SIMPLE methods.
While for the primal system, we employ a tailored multigrid method specifically developed for the unfitted meshes.
Lastly, the comparison between the proposed preconditioners is made through several numerical experiments.

 \end{abstract}

\section{Introduction}
In the modeling of many real-world engineering problems, we encounter material discontinuities, such as inclusions.
The inclusions are found naturally in the materials, or can be artificially introduced to produce desired mechanical behavior.
The finite element (FE) modeling of such inclusions requires to generate meshes that can resolve the interface between the matrix and inclusions, which can be a computationally cumbersome and expensive task.
In order to avoid the need of generating such fitted meshes, the extended finite element method (XFEM) was introduced~\cite{sukumar_modeling_2001}.
The XFEM method can be categorized as an unfitted finite element method.
For the unfitted FE methods, it is not necessary to have the meshes that resolve the interfaces exactly.
In this framework, structured uniform meshes are typically used and the associated FE spaces are enriched to capture the interface information.
There are multiple ways for enriching the FE spaces and for enforcing the constraints on the interfaces.
Here, we employ the method of Lagrange multipliers to enforce the interface conditions.
The method of Lagrange multipliers is quite robust, but it gives rise to the mixed formulation and requires the solution of the saddle point type linear system of equations.

In this work, we model a domain with multiple inclusions using a variant of the unfitted FE method.
We briefly introduce the unfitted discretization methods and then discuss the solution strategies for solving the arising saddle-point systems.
In particular, we present a survey on the available preconditioning strategies for solving the primal and dual systems.
 
\section{Model problem}
We assume a domain $\Omega \in \bR^2$ with the Lipschitz continuous boundary $\partial \Omega$.
This domain $\Omega$ is assumed to be decomposed into distinct non-overlapping parts, the matrix $\Omega_M$ and $n$ inclusions $\{\Omega_I^1,\Omega_I^2,\ldots,\Omega_I^n\}$.
The inclusions are assumed to be completely embedded in the matrix domain $\Omega_M$.
Each inclusion is separated from the domain $\Omega_M$ by interfaces $\{\Gamma^1,\Gamma^2,\ldots,\Gamma^n\}$.
The interfaces associated with the inclusions are assumed to be sufficiently smooth.
Thus, the domain is defined as, $\Omega = \bigcup_i^n (\Omega_I^i \cup \Gamma^i) \cup \Omega_M$.
The boundary $\partial \Omega$ is decomposed into Dirichlet and Neumann boundaries given as $\partial \Omega_D$ and $\partial \Omega_N$, respectively.
We assume that the domain is subjected to volume forces $\vecf$ and the traction/surface forces ${\vect_N}$.
Under the influence of these external forces, the domain under goes deformation, denoted as $\vecu:=\{\vecu_M,\vecu_I^1,\vecu_I^2,\ldots,\vecu_I^n\}$.
The deformation or the displacement field is defined as a sufficiently regular function for the matrix domain as $\vecu_M$ and the inclusions as $\vecu_I^i$.
On the interfaces, we assume the perfect bonding condition, such that the displacement field is defined to be continuous.
The equation of the equilibrium, the boundary conditions and the interface conditions on such domain are written as
\begin{equation}
  \begin{aligned}
    \nabla \cdot \bs{\sigma} + \vecf & = 0       &  & \text{in } \Omega,                                    \\
    \vecu                            & = 0       &  & \text{on } \partial \Omega_D,                         \\
    \bs{\sigma} \vecn                & = \vect_N &  & \text{on } \partial \Omega_N,                         \\
    \llbracket \vecu_i \rrbracket    & = 0       &  & \text{on } \Gamma^i, \text{ for } i=\{1,2,\ldots,n\},
  \end{aligned}
  \label{eq:linear_elasticity}
\end{equation}
where, $\bs{\sigma}$ denotes the Cauchy stress tensor and $\vecn$ denotes the outward normal.
The jump of the displacement field is defined as $\llbracket \vecu_i \rrbracket := \vecu_M \vert_{\Gamma_i} - \vecu_I^i\vert_{\Gamma_i}$.
We assume that the material is linear elastic and the constitutive law is provided by Hooke's law
\(
\bs{\sigma} = \lambda \text{tr}(\bs{\varepsilon}) \matI + 2 \mu \bs{\varepsilon},
\label{eq:contitutive_law}
\)
where $\lambda$ and $\mu$ are Lam\'e parameters, and $\text{tr}(\cdot)$ denotes the trace operator.
In this work, we consider the linearized strain tensor $\bs{\varepsilon}$, given as $\bs{\varepsilon} := \frac{1}{2}\big( \nabla \vecu +(\nabla \vecu)^\T \big)$.
 
\subsection{Unfitted discretization}
\begin{figure*}[t]
  \centering
  \subcaptionbox{Inclusions in a matrix}[.24\textwidth]{\includegraphics{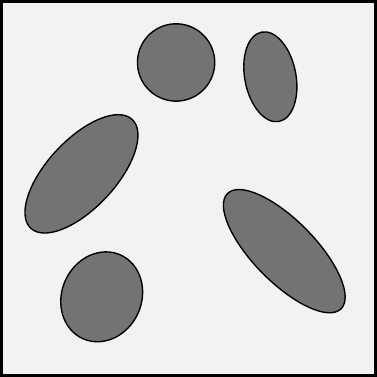}
  }
  \hfill \subcaptionbox{An unfitted mesh $\widetilde{\pT}_h$}[.24\textwidth]{\includegraphics{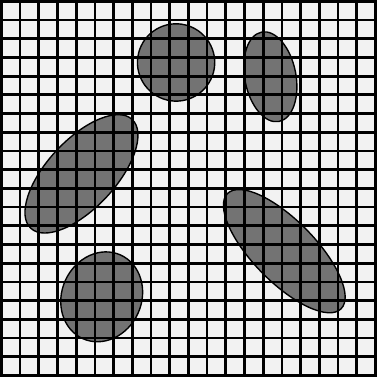}
  }
  \hfill
  \subcaptionbox{$\pT_{h,M}$ encapsulating the matrix $\Omega_M$}[.24\textwidth]{\includegraphics{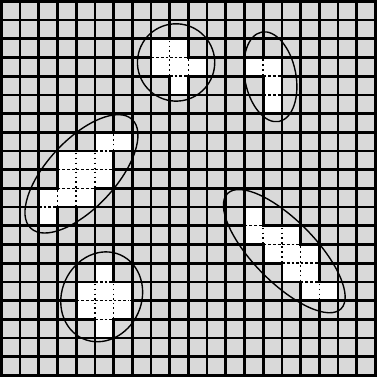}
  }
  \hfill
  \subcaptionbox{$\pT^i_{h,I}$ encapsulating the inclusions $\Omega^i_I$}[.24\textwidth]{
  \includegraphics{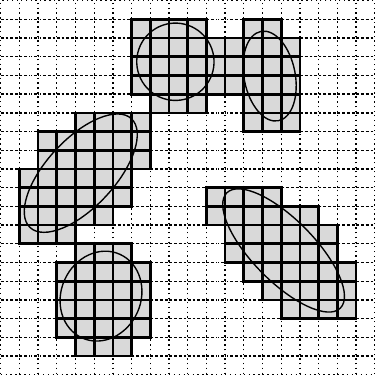}
    }
  \caption[2D mesh]{A domain with multiple inclusions and the unfitted meshes which encapsulates the matrix and inclusions.}
  \label{fig:domain_decomp}
\end{figure*}

We use a Cartesian mesh $\widetilde{\pT}_h$, which is assumed to be quasi-uniform and shape regular.
The mesh $\widetilde{\pT}_h$ is fitted to the domain $\Omega$, but it is not fitted with the interfaces $\Gamma^i$.
We assume, the interfaces $\Gamma^i$ are resolved sufficiently well by the mesh $\widetilde{\pT}_h$ and the curvature of the interfaces is bounded.
The interfaces intersect with a given edge only once and do not pass through the nodes.
The mesh $\widetilde{\pT}_h$ is treated as a background mesh that captures the matrix domain and all inclusions.
Now, we can define the so-called active mesh which is associated with either the matrix subdomain or inclusions.
The active mesh is strictly intersected by the subdomains $\Omega_M$ or $\Omega_I^i$, given as
\[
  \pT_{h,M} = \{ K\in \widetilde{\pT}_h: K\cap \Omega_M  \neq \emptyset\}, \quad  \pT_{h,I}^i = \{ K\in \widetilde{\pT}_h: K\cap \Omega_I^i  \neq \emptyset\}, \text{ for } i=\{1,2,\ldots,n\}.
\]
Thus, each subdomain is encapsulated by the respective active mesh, such that all the elements that do not intersect with the subdomain or the interfaces are excluded.
In Figure~\ref{fig:domain_decomp}, we can see an example of active meshes associated with the matrix and inclusions.
The set of elements that are intersected by the interfaces are defined as,
\[
  \pT_{h,\Gamma}^i = \{K \in \widetilde{\pT}_h: K\cap \Gamma^i \neq \emptyset \} \quad \text{ for } i=\{1,2,\ldots,n\}.
\]
The interface meshes $\pT_{h,\Gamma}^i$ are doubled and are associated with the respective matrix and inclusion subdomains, which we denote by $\pT_{h,\Gamma_I}^i$ and $\pT_{h,\Gamma_M}^i$, respectively.
Thus, we enrich the mesh with the extra degrees of freedom (DOFs) in the vicinity of the interfaces.

We define a finite element (FE) space over the background mesh $\widetilde{\pT}_h$ as
\[
  \label{eq:background_FEspace}
  \widetilde{\bpV}_h = \{\vecv \in \matH^1 (\widetilde{\pT}_h): \vecv\vert_K \in \pQ_1(K), \vecv\vert_{(\partial \widetilde{\pT}_h)_D} = 0, \forall K\in \widetilde{\pT}_h \},
\]
where $\pQ_1$ denotes the space of piecewise bilinear functions.
We define characteristic functions on each computational subdomains $\Omega_M$, $\chi_{\Omega_M}: \bR^d \to \bR$ and $\Omega_I^i$, $\chi_{\Omega_I^i} : \bR^d \to \bR$ as
\[
  \label{eq:Heaviside}
  \chi_{\Omega_M}(\matX) = \begin{cases}
    1 & \quad \forall \matX \in \overline{\Omega_M}, \\
    0 & \quad \text{otherwise,}
  \end{cases} \text{ and }\
  \chi_{\Omega_I^i}(\matX) = \begin{cases}
    1 & \quad \forall \matX \in \overline{\Omega_I^i}, \\
    0 & \quad \text{otherwise,}
  \end{cases}
  \quad \text{ for } i=\{1,2,\ldots,n\}.
\]
These characteristic functions are employed to restrict the support of the finite element space $\widetilde{\bpV}_h$ to the respective subdomains $\bpV_{h,M} = \chi_{\Omega_M}(\matX)\widetilde{\bpV}_h$ and $\bpV_{h,I}^i = \chi_{\Omega_I^i}(\matX)\widetilde{\bpV}_h$.

Now, we seek the solution $\vecu_h = (\bigoplus_{i=1}^n \vecu_{h,I}^i \oplus \vecu_{h,M})$ in FE space $\bpV_h = (\bigoplus_{i=1}^n\ \bpV_{h,I}^i \oplus \bpV_{h,M})$.
The enriched FE space $\bpV_h$ is equipped with the extra DOFs and it also accommodates the jumps across the interfaces as FE space associated with each subdomain is decoupled.
We also note that the nodal basis functions spanning the enriched FE space also have restricted support, which means the support of the basis functions vanishes outside of the subdomain.

\subsection{Variational formulation}
In order to derive the variational formulation of the problem \eqref{eq:linear_elasticity}, we utilize the enriched FE space $\bpV_h$ and the method of Lagrange multipliers to enforce the interface condition.
We define the multiplier space $\bpM_h$ as the dual of the trace space of $\bpV_h$, thus $\bpM_h \subset \matH^{-\half}(\Gamma)$.
The variational formulation of the problem \eqref{eq:linear_elasticity} is given as, find $(\vecu_h,\bs{\lambda}_h) \in \bpV_h\times\bpM_h$ such that
\begin{equation}
  \begin{aligned}
     & a(\vecu_h,\vecv_h) + b(\vecv_h,\bs{\lambda}_h) & = & \ F(\vecv_h) \quad & \forall \vecv_h \in \bpV_h,    \\
     & b(\vecu_h,\bs{\mu}_h)                          & = & \ 0          \quad & \forall \bs{\mu}_h \in \bpM_h.
    \label{eq:weak_form}
  \end{aligned}
\end{equation}
The bilinear and the linear forms are defined as
\[
  \begin{aligned}
    \label{eq:elast_bilinear_linear}
    a(\vecu_h,\vecv_h) & := \int_{\Omega_M} \bs{\sigma}(\vecu_{h,M}):\bs{\varepsilon}({\vecv_h}) d\Omega+ \sum_{i=1}^n \int_{\Omega_I^i} \bs{\sigma}(\vecu_{h,I}^i):\bs{\varepsilon}({\vecv_h}) d\Omega,                                                   \\
    F(\vecv_h)         & := \int_\Omega \bs{f} \cdot \vecv_h d\Omega + \int_{\partial \Omega_N}  \vect_N \cdot \vecv_h d\Gamma, \quad b(\vecu_h,\bs{\bs{\mu}_h}) := \sum_{i=1}^n \int_{\Gamma_i} \llbracket \vecu_{h,i} \rrbracket\ \bs{\mu}_h \ d \Gamma,
  \end{aligned}
\]
where $a(\cdot,\cdot): \bpV_h \times \bpV_h \to \bR$ is as symmetric continuous coercive bilinear form, ${F(\cdot):\bpV_h \to \bR}$ denotes continuous linear form, and $b(\cdot,\cdot):\bpV_h \times\bpM_h \to \bR$ is a bilinear form.
In the next section, we discuss the Lagrange multiplier space for the unfitted FEM framework.

\subsubsection{A stable Lagrange multiplier space}
The weak formulation given in \eqref{eq:weak_form} is stable  only if the following discrete inf-sup condition is satisfied,
\[
  \inf_{\bs{\lambda}_h \bpM_h} \sup_{\vecu_h \in \bpV_h} \frac{b(\vecu_h,\bs{\lambda}_h)}{\|\bs{\lambda}_h\|_{\bpM_h} \|\vecu_h\|_{\bpV_h} } \geqslant \beta > 0,
\]
where the constant $\beta$ is independent of the mesh size $h$.
Thus, the choice of the multiplier space is quite essential for ensuring the stability of the discretization method.

In this work, we construct the multiplier space by employing the vital vertex algorithm~\cite{bechet_stable_2009}.
In the unfitted FE framework, if the multiplier space is constructed by using all nodes that are associated with the interface mesh $\pT_{h,\Gamma}^i$, we obtain a very rich multiplier space that does not satisfy the inf-sup condition.
Thus, the naive method for constructing the multiplier space leads to an unstable discretization method.
The vital vertex algorithm selects a set of vertices from all the vertices in such a way that the resultant multiplier space is not too rich.
In addition, the basis functions associated with the multiplier space are constructed as trace of the FE basis functions from the primal space and we not require to construct the lower dimensional basis functions for numerical integration on  interfaces.

\subsubsection{Ghost penalty stabilization}
In the unfitted methods, the background mesh and the interfaces are allowed intersect arbitrarily.
Due to this flexibility, in some case the elements are cut into disproportionate fractions, which could affect aversely the condition number of the system matrices.
The ghost penalty stabilization was introduced in order to alleviate the issue of ill-conditioning, and regain the control of the gradients over the cut elements with very small support~\cite{burman_ghost_2010}.
We define the set of faces $\pG_{h,\Gamma}$ for each subdomain associated with the inclusions $\Omega_I^i$ and the matrix subdomain $\Omega_M$, as
\[
  \begin{aligned}
    \label{eq:}
    \pG^i_{h,\Gamma_I} & =\{ G \subset \partial K \mid K\in (\pT^i_{h,\Gamma_I}),\ \partial K \cap \partial \pT_I^i = \emptyset\} \quad \text{for } i=\{1,2,\ldots,n\} , \\
    \pG^i_{h,\Gamma_M} & =\{ G \subset \partial K \mid K\in (\pT^i_{h,\Gamma_M}),\ \partial K \cap \partial \pT_M = \emptyset\} \quad \text{for } i=\{1,2,\ldots,n\}.
  \end{aligned}
\]
The ghost penalty term is enforced on the set of edges $\pG_{h,\Gamma}$, and it is defined as
\[
  j(\vecu_h,\vecv_h) = \sum_{D\in\{M,I\}} \sum_{i=1}^n \sum_{G\in \pG^i_{h,\Gamma_{D}}} \int_G \epsilon_G \frac{h_G}{(2\mu_{D}+\lambda_{D})} \llbracket \nabla_{\vecn_G}\pE_{h,D}\vecu_h \rrbracket \llbracket \nabla_{\vecn_G}\pE_{h,{D}} \vecv_h \rrbracket \ dG,
\]
where $h_G$ is the diameter of face $G$, $\vecn_G$ denotes unit normal to face $G$, $\epsilon_G$ is a positive constant, and $\lambda_D, \mu_D$ denotes Lam\'e parameters associated with either inclusions or the matrix subdomains.
Here, $\pE_{h,D}$ denotes the canonical extension of the function from the domain to the background mesh, which is defined as $\pE_{h,{D}}:\bpV_h\vert_{K_{D}} \to \widetilde{\bpV}_h\vert_K$.
By adding the ghost penalty term to the bilinear form $a(\cdot,\cdot)$, we can regain the control over the gradients for the small cut elements.
The modified or stabilized weak form for our model problem can be given as:
\begin{equation}
  \begin{aligned}
     & a_j(\vecu_h,\vecv_h) + b(\vecv_h,\bs{\lambda}_h) &  & = F(\vecv_h) \quad &  & \forall \vecv_h \in \bpV_h,    \\
     & b(\vecu_h,\bs{\mu}_h)                            &  & = 0 \qquad         &  & \forall \bs{\mu}_h \in \bpM_h,
    \label{eq:weak_form_2}
  \end{aligned}
\end{equation}
where the bilinear form $a_j(\cdot,\cdot)$ is defined as $ a_j(\vecu_h,\vecv_h) = a(\vecu_h,\vecv_h) + j(\vecu_h,\vecv_h)$.
 
\section{Solution strategy}
In this section, we discuss the solution strategies for solving the saddle point formulation \eqref{eq:weak_form_2}.
The saddle point problem in the algebraic form is given as,
\begin{equation}
  \begin{aligned}
    \label{eq:spp_matvec}
    \matA \vecx + \matB^\T \vecy & = \vecb   \\
    \matB \vecx                  & = \bs{0},
  \end{aligned}
\end{equation}
where $\matA \in \bR^{n\times n}$, $\vecb \in \bR^n$, $\matB \in \bR^{m \times n}$, $\text{rank}(\matB)=m$, $m \ll n$, with the unknowns $\vecx \in \bR^n$ and the Lagrange multipliers $\vecy \in \bR^m$.
The matrix $\matA$ is a symmetric positive semidefinite, and it does not have an exact inverse due the floating domains associated with the inclusions.
The problem \eqref{eq:spp_matvec} is solvable if and only if, ${(\vecb - \matB^\T \vecy) \perp \Ker(\matA)}$, where $\Ker(\matA)$ denotes the kernel of the matrix $\matA$.
In practice, special care is required in order to solve the system with a non-trivial kernel.
In particular, if the null space of the matrix $\matA$ is known, an iterative process can create a sequence of iterates that are orthogonal to $\Ker(\matA)$~\cite{bochev_finite_2005}.
In this work, we take a different approach and transform the block diagonal matrix $\matA$ into an equivalent symmetric positive definite matrix.
To this end, we use the method of augmented Lagrangian to reformulate the problem \eqref{eq:spp_matvec} and replace the problem with an equivalent saddle-point formulation,
\begin{equation}
  \label{eq:saddle_AL}
  \begin{aligned}
    \matA_\gamma \vecx + \matB^\T \vecy & = \vecb   \\
    \matB \vecx                         & = \bs{0},
  \end{aligned}
\end{equation}
where $\matA_\gamma := \matA + \gamma_s \matB^\T \matB $ is the augmented matrix and $\gamma_s$ denotes a stabilization parameter.
Now, the matrix $\matA_\gamma$ in the saddle point problem \eqref{eq:saddle_AL} is a positive definite matrix.
We choose the value of the penalty parameter as
\(
\gamma_s := \frac{\lambda_{\max} (\matA)}{\lambda_{\max} (\matB^\T \matB)},
\)
where $\lambda_{\max}$ denotes the largest eigenvalue of the corresponding matrix.
This penalty parameter is quite an attractive option as it minimizes the condition number of the whole saddle point system and by extension enhances the convergence of the standard iterative methods~\cite{greif_augmented_2004}.
From now onwards, all the methods discussed in this section are concerned with solving the modified saddle point problem \eqref{eq:saddle_AL}.

\subsection{Uzawa methods}
The Uzawa method is widely used for solving saddle point problems arising from the Stokes problem, incompressible solid and fluid mechanics problems.
It consists of the following coupled iteration
\begin{equation}
  \begin{aligned}
    \label{eq:uzawa}
    \matA_\gamma \vecx^{(k+1)} & = \vecf - \matB^\T \vecy^{(k)},               \\
    \vecy^{(k+1)}              & = \vecy^{(k)} + \omega (\matB \vecx^{(k+1)}),
  \end{aligned}
\end{equation}
where $\omega > 0$ is a relaxation parameter.
Here, the first equation can be used to eliminate $\vecx^{(k+1)}$ from the second equation, which leads to a stationary Richardson iterative method for the dual system
\(
\vecy^{(k+1)} = \vecy^{(k)} + \omega (\matB \matA_\gamma^{-1} ( \vecf - \matB^\T \vecy^{(k)} )).
\)
The Uzawa methods are generally slow to converge, as their convergence rate depends on the choice of the relaxation parameter $\omega$.
Therefore, it is desirable to replace the Richardson method in the Uzawa scheme with more robust steepest descent or conjugate gradient (CG) method~\cite{dietrichbraess2007-04-28}.
Simultaneously, we also employ some preconditioning strategies to accelerate the convergence of the Uzawa-CG method and to reduce the overall computational cost of the solution scheme \eqref{eq:uzawa}.
The preconditioned Uzawa method can be written as
\begin{equation}
  \begin{aligned}
    \label{eq:precond_uzawa}
    \matA_\gamma \vecx^{(k+1)} & = \vecf - \matB^\T \vecy^{(k)}                  \\
    \vecy^{(k+1)}              & = \vecy^{(k)} + \matP^{-1} \matB \vecx^{(k+1)},
  \end{aligned}
\end{equation}
where $\matP \in \bR^{m \times m}$ denotes the preconditioner matrix for the dual system.

\subsection{Preconditioners for dual systems}
In this section, we describe some of the preconditioning techniques for dual systems.
Since the saddle point problem appears in many practical applications, there have been many efforts for developing optimal preconditioning strategies.
In this section, we focus on the preconditioners developed in two different areas, namely SIMPLE and FETI methods.

\subsubsection{SIMPLE Preconditioner}
Semi-Implicit Method for Pressure Linked Equation (SIMPLE) is a solution technique developed for the problems arising from Navier-Stokes equations~\cite{s.v.patankar1979-12-31}.
The incompressible Navier-Stokes problems also gives rise to the saddle point type linear system of equations as~\eqref{eq:precond_uzawa}.
In the SIMPLE method, $\vecx$ represents the velocity vector and $\vecy$ represents a pressure vector, while the matrix $\matB$ represents the constraint on the pressure vector.
The original formulation of the SIMPLE method can be viewed as a semi-implicit Uzawa method and the preconditioner is normally used for a coupled iteration.
Due to the structural similarity of the problems, in this work, we aim to employ the SIMPLE preconditioner for solving the dual problem.
The SIMPLE preconditioner is given as
\begin{equation}
  \label{eq:simple_preconditioner}
  \matP_{S}^{-1}:=(\matB\matD^{-1}\matB^\T)^{-1}.
\end{equation}

\subsubsection{FETI Preconditioners}
The Finite Element Tearing and Interconnecting (FETI) methods are introduced by Farhat and Raux as a non-overlapping domain decomposition methods~\cite{farhat_method_1991}.
FETI is a group of iterative sub-structuring methods for solving large systems of linear equations arising from FEM discretization.
By design, the FETI methods are parallel solution methods, where the computational domain is decomposed into multiple subdomains.
These subdomains are distributed among multiple processors.
On each processor, a local Neumann problem is solved along with a coarse problem, used for global information transfer.
The continuity between all subdomains is imposed by means of the method of Lagrange multipliers.
We aim to leverage some preconditioners developed for the FETI methods, as the algebraic formulation of the FETI method also gives rise to the saddle point system \eqref{eq:saddle_AL}.

In earlier FETI literature~\cite{farhat_optimal_1994}, Farhat et al.~proposed the Dirichlet preconditioner, given as
\begin{equation}
  \label{eq:feti_d_preconditioner}
  \matP_{D}^{-1}:=\matB \matA \matB^\T.
\end{equation}
This preconditioner, even though being useful, is quite elementary.
More complex and better alternatives for the preconditioners are proposed by Lacour~\cite{lacour_iterative_1996}, given as
\begin{equation}
  \label{eq:feti_preconditioner1}
\matP_{L}^{-1}:= (\matB\matB^\T)^{-1}\matB \matA \matB^\T (\matB\matB^\T)^{-1},
\end{equation}
and by Klawonn and Widlund \cite{klawonn_feti_2001}, given as
\begin{equation}
  \label{eq:feti_preconditioner2}
  \matP_{F}^{-1}:=(\matB\matD^{-1}\matB^\T)^{-1}\matB \matD^{-1} \matA \matD^{-1} \matB^\T (\matB\matD^{-1}\matB^\T)^{-1},
\end{equation}
where, $\matD = \text{diag}(\matA)$.

The SIMPLE and FETI preconditioners are very attractive possibilities for solving the dual systems arising from the saddle point system~\eqref{eq:spp_matvec}.
We aim to use these preconditioners for solving the dual problems arising in the Uzawa method \eqref{eq:precond_uzawa}.

\subsection{Multigrid method for the primal system}
For solving the primal problem, we rely on the multigrid methods.
The multigrid methods are ideal preconditioners for solving large systems arising from the discretization of partial differential equations as they have optimal complexity.
The multigrid method can be expressed as an ideal combination of the smoothing iterations and the coarse level corrections.
In the fitted FE methods, a hierarchy of nested meshes is generated and the standard interpolation operator and its adjoint are used as transfer operators.
In the unfitted methods, generally, the structured meshes are used as background meshes, and a hierarchy of nested background meshes are generated.
Unfortunately, we can not use the standard interpolation operators as the transfer operator in multigrid method for the unfitted FE methods.
This is due to the fact that, the locations of the interfaces are not fixed with respect to the background meshes, and encapsulated meshes might not be nested even though the background meshes are nested.
Thus, a hierarchy of non-nested meshes in turn gives rise to a hierarchy of non-nested FE spaces.
In order to create a hierarchy of nested FE spaces from the hierarchy of non-nested meshes, we employ the $L^2$-projection based transfer operators.

We propose to employ tailored pseudo-$L^2$-projection-based transfer operators for the multigrid method for unfitted discretization method.
It has been shown in~\cite{kothari_multigrid_2019}, the multigrid methods equipped with such transfer operator are robust for solving the linear system arising in unfitted finite element methods.
In this work, we use this robust multigrid method as a preconditioner for the CG method.

\section{Numerical results}
In this section, we evaluate the performance of the discretization method and the preconditioners for the primal and dual systems.

\begin{figure*}[]
  \centering
  \subcaptionbox{Problem setup\label{subfig:problemsetup}}[.24\textwidth]{\includegraphics{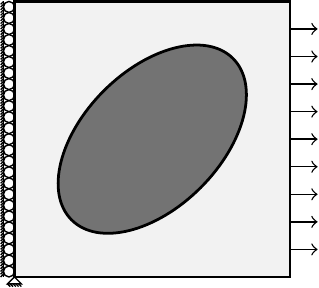}
    }
  \hfill \subcaptionbox{Soft inclusion}[.24\textwidth]{\includegraphics[width=0.25\textwidth]{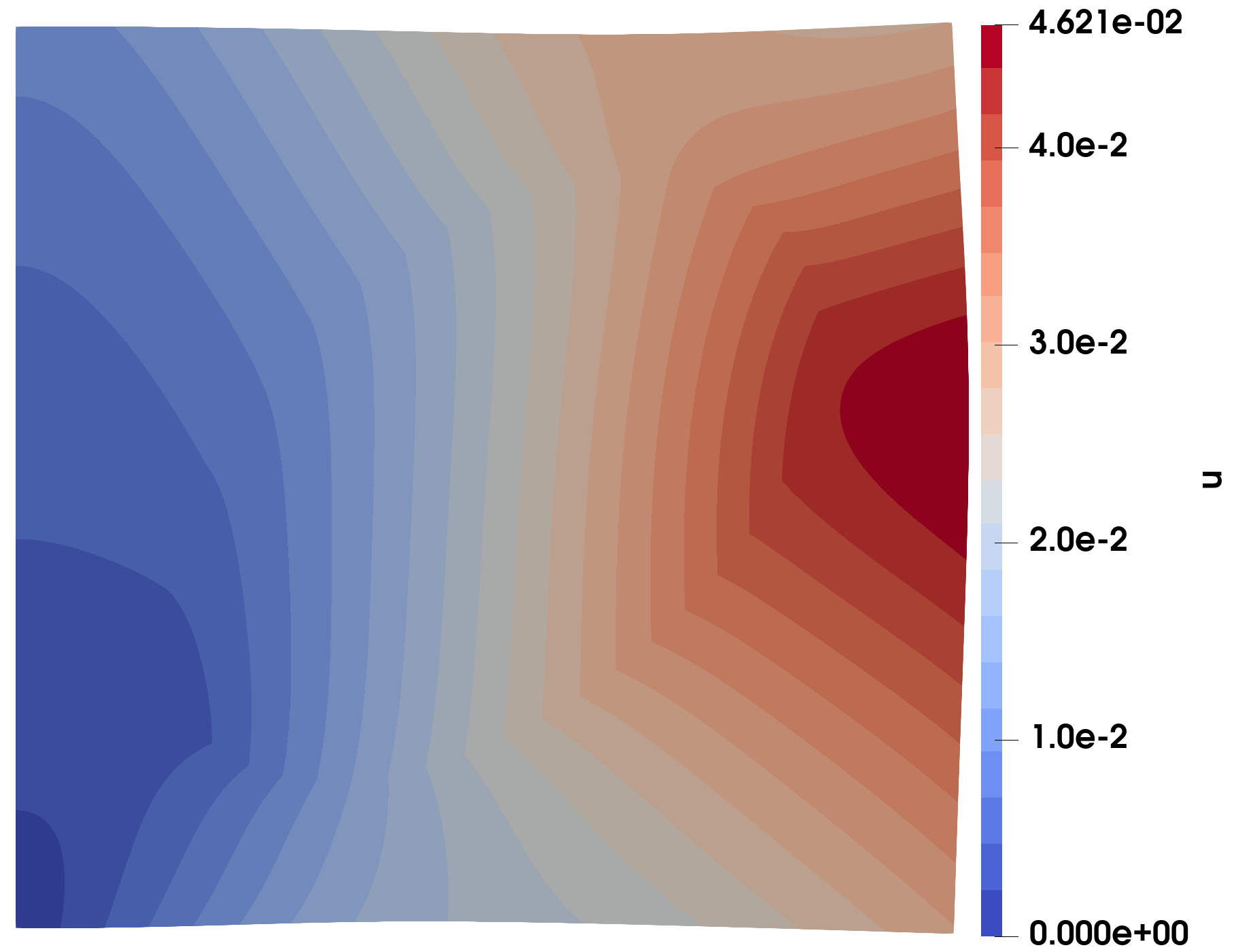}
  }
  \hfill
  \subcaptionbox{Hard inclusion}[.24\textwidth]{\includegraphics[width=0.25\textwidth]{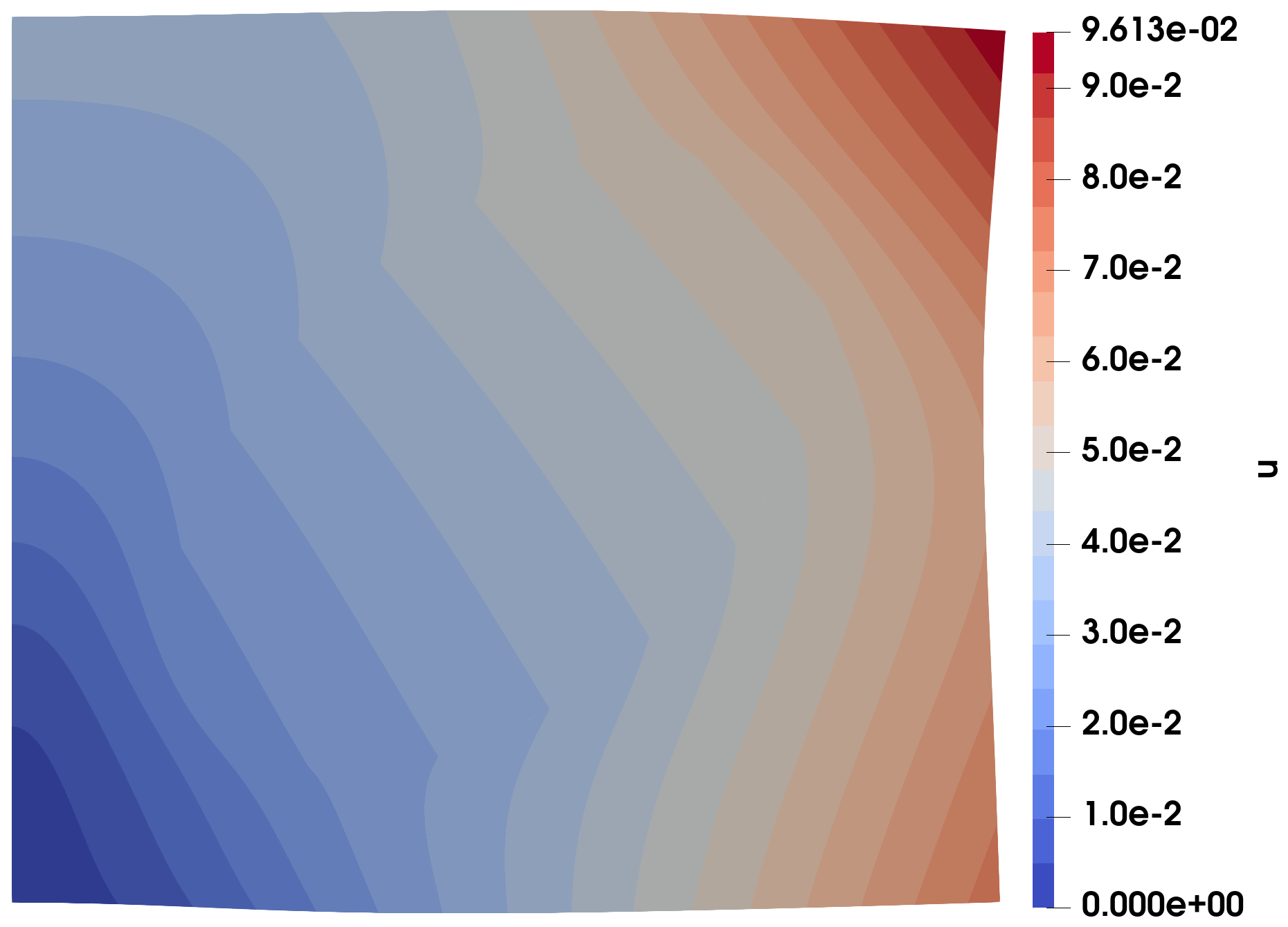}
  }
  \hfill
  \subcaptionbox{Uniform inclusion}[.24\textwidth]{\includegraphics[width=0.25\textwidth]{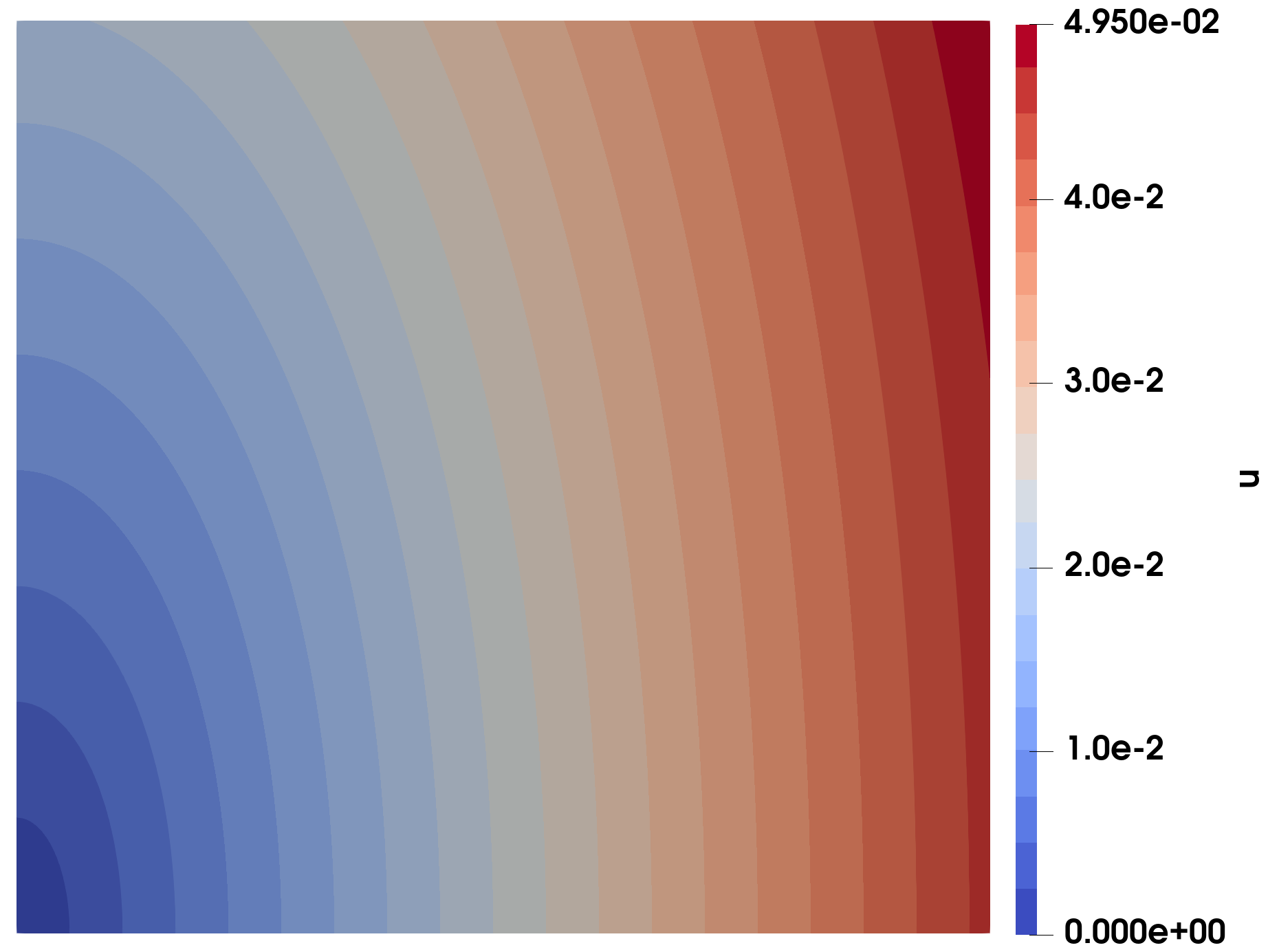}
  }
  \caption[]{The problem setup and the resultant displacement field for different test cases.}
  \label{fig:resuting_displacement}
\end{figure*}

We consider a domain $\Omega = [0,1]^2$, with an elliptical inclusion $\Omega_I$, subjected to Dirichlet boundary condition on the left edge and Neumann boundary condition on the right edge.
The magnitude of the traction on the right edge is given as $10^4$.
The configuration of the problem is shown in Figure~\ref{subfig:problemsetup}.
The elliptical interface is defined as a zero level set of function $\Lambda(\hat{x},\hat{y}):=1 - (\hat{x}^2/a^2) - (\hat{y}^2/b^2)$, where $a,b$ denote the semi-major and semi-minor axes, and the rotated coordinates with orientation $\theta$ are given as $\hat{x}=(x-c_x)\cos(\theta) + (y -c_y)\sin(\theta)$ and $\hat{y}=-(x-c_x)\sin(\theta) + (y-c_y)\cos(\theta)$.
Here the center of the ellipse $(c_x,c_y)$ are given as $(0.5,0.5)$, the orientation is chosen as $\theta=\pi/4$ and semi-major and semi-minor axes are given as $a = r, b = 0.6r$, where $r=\sqrt{3-2\sqrt{2}}$.
We define different Young's modulus for the inclusion and the matrix, and set up three different test cases: Soft inclusion ($E_M = {4\times10^5}\,\text{MPa}, E_I =  {10^5}\,\text{MPa}$), Hard inclusion ($E_M = 10^5\,\text{MPa}, E_I =  {4\times10^5}\,\text{MPa}$), and Uniform inclusion ($E_M = E_I = {2\times10^5}\,\text{MPa}$), while the Poisson ratio is chosen as $\nu=0.3$ for all test cases.
The resultant displacement field for all three test cases can be seen in Figure~\ref{fig:resuting_displacement}.

\subsection{Study of the discretization error}
We start with a mesh that has 50 elements in each direction, which we denote as $L1$.
We uniformly refine the mesh $L1$ to obtain a background mesh hierarchy denoted as $L2,\ldots, L5$.
This mesh hierarchy is used for measuring the discretization error and as a multilevel hierarchy in the multigrid method.

We compare the performance of the unfitted FE discretization by investigating the convergence of the numerical approximation.
The numerical error is defined as $\vece_h = \vecu_{f}-\vecu_h$, where the subscript $f$ denotes the approximation on the finest level $L5$, and $h$ denotes the approximation on the current level.
The error is computed in the $L^2$-norm, $H^1$-seminorm and energy norm, which are defined as
$ \| \vece_h \|^2_{L^2(\Omega)} := \sum_{K\in \pT_h} \int_K ( \vece_h )^2 \,d\Omega$,
$ \vert \vece_h \vert^2_{H^1(\Omega)} := \sum_{K \in \pT_h} \int_K ( \nabla \vecu_{f}- \nabla \vecu_h )^2 \,d\Omega$,
$ \| \vece_h \|^2_E :=  \sum_{K\in \pT_h} \int_K \bs{\sigma}(\vece_h ):\bs{\varepsilon}(\vece_h ) \,d\Omega$.
From the Figure \ref{fig:convergence_1}, we can see that for all test cases, the rate of convergence of the discretization error in $L^2$-norm is $\pO(h^2)$, for the $H^1$-seminorm and energy norm is $\pO(h)$.

\tikzset{external/export next=false}
\begin{figure}[hbt!]
  \centering
  \begin{tikzpicture}[]
    \begin{groupplot}[
        group style={
group size = 3 by 1,
x descriptions at=edge bottom,
horizontal sep=0.6cm,
},
scale=0.54,xlabel=$\log_{10}(h)$, minor tick num=4, grid=major,tick style={thick},grid style={thick,dashed,gray},
label style={font=\scriptsize}, tick label style={font=\scriptsize}, legend style={font=\scriptsize}
      ]
\nextgroupplot[align=left, title={{\scriptsize Soft inclusion}},ylabel={$\log_{10}(\text{Error})$}]
      \addplot [very thick,red!70!black,mark=diamond*,every mark/.append style={solid}]table[x=h, y=l2_soft, col sep=comma]{example1.csv};
      \label{pgfplots:1L2}
      \addplot [densely dashed,very thick,green!70!black,mark=triangle*,every mark/.append style={solid}]table[x=h, y=h1_semi_soft, col sep=comma]{example1.csv};
      \label{pgfplots:1SH1}
      \addplot [densely dotted,very thick,blue!70!black,mark=*,every mark/.append style={solid}]table[x=h, y=energy_soft, col sep=comma]{example1.csv};
      \label{pgfplots:1Energy}

      \nextgroupplot[align=left, title={{\scriptsize Hard inclusion}},ylabel={}]
      \addplot [very thick,red!70!black,mark=diamond*,every mark/.append style={solid}]table[x=h, y=l2_hard, col sep=comma]{example1.csv};
\addplot [densely dashed,very thick,green!70!black,mark=triangle*,every mark/.append style={solid}]ttable[x=h, y=h1_semi_hard, col sep=comma]{example1.csv};
\addplot [densely dotted,very thick,blue!70!black,mark=*,every mark/.append style={solid}]table[x=h, y=energy_hard, col sep=comma]{example1.csv};

      \nextgroupplot[align=left, title={{\scriptsize Uniform inclusion}},ylabel={}]
      \addplot [very thick,red!70!black,mark=diamond*,every mark/.append style={solid}]table[x=h, y=l2_equal, col sep=comma]{example1.csv};
\addplot [densely dashed,very thick,green!70!black,mark=triangle*,every mark/.append style={solid}]table[x=h, y=h1_semi_equal, col sep=comma]{example1.csv};
\addplot [densely dotted,very thick,blue!70!black,mark=*,every mark/.append style={solid}]table[x=h, y=energy_equal, col sep=comma]{example1.csv};

    \end{groupplot}
    \matrix [ draw, matrix of nodes, anchor = west, node font=\scriptsize,
    column 1/.style={nodes={align=center,text width=1.2cm}},
    column 2/.style={nodes={align=center,text width=0.6cm}},
    ] at ($(group c3r1)+(2,0)$)
    {
    $\|\vece_h\|_{L^2(\Omega)}$         & \ref{pgfplots:1L2}  \\
    $\vert \vece_h \vert_{H^1(\Omega)}$  & \ref{pgfplots:1SH1} \\
    $\|\vece_h\|_{E}$         & \ref{pgfplots:1Energy} \\
    };
\end{tikzpicture}
  \caption{Convergence of error in $L^2$-norm, $H^1$-seminorm and energy norm.}
  \label{fig:convergence_1}
\end{figure}
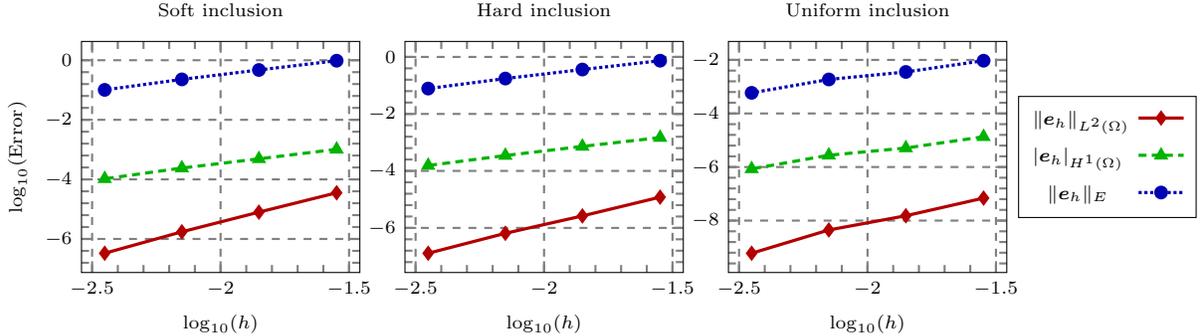
 
\subsection{Comparison of preconditioners}
In this section, we compare the performance of the various preconditioners for the primal and dual problems in the saddle-point system.
First, we compare the condition numbers of the preconditioned Schur-complement matrix ($\matS := \matB \matA^{-1}\matB^\T$) with respect to decreasing the mesh size.
The result of this comparison is shown in Figure~\ref{fig:convergence_1_cond}, where we can see that for all test cases, the preconditioner $\matP_F$ is the most stable amongst all other options.
The Dirichlet preconditioner causes the condition number of the preconditioned system to increase rather than decrease, and the SIMPLE preconditioner is effective only for the uniform inclusion case.
The most effective preconditioners are definitely $\matP_L$ and $\matP_F$, as the condition number of the preconditioned system is extremely stable and it does not increase with the decreasing mesh size.

\tikzset{external/export next=false}
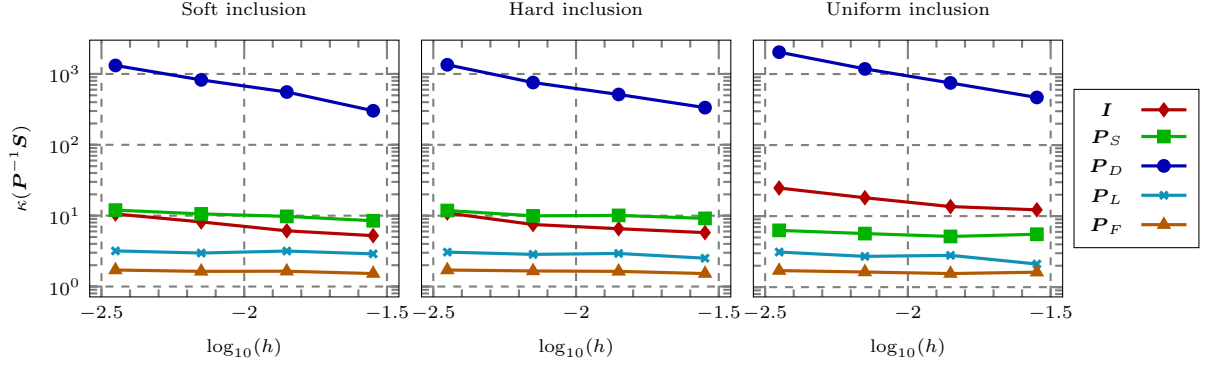
\begin{figure}[t]
  \centering
  \begin{tikzpicture}[]
    \begin{groupplot}[
        group style={
group size = 3 by 1,
x descriptions at=edge bottom,
horizontal sep=0.3cm,
},
scale=0.6, xlabel=$\log_{10}(h)$,ymode=log, minor tick num=4, grid=major,tick style={thick},grid style={thick,dashed,gray},
        ymin=0, ymax=3e3,
        label style={font=\scriptsize}, tick label style={font=\scriptsize}, legend style={font=\scriptsize}
      ]
\nextgroupplot[align=left, title={{\scriptsize Soft inclusion}},ylabel={$\kappa(\matP^{-1}\matS)$}]
      \addplot [very thick,red!70!black,mark=diamond*,every mark/.append style={solid}]table[x=h, y=s_soft, col sep=comma]{example1_cond.csv};
      \label{pgfplots:s_soft}
      \addplot [very thick,green!70!black,mark=square*,every mark/.append style={solid}]table[x=h, y=simple_soft, col sep=comma]{example1_cond.csv};
      \label{pgfplots:simple_soft}
      \addplot [very thick,blue!70!black,mark=*,every mark/.append style={solid}]table[x=h, y=dirichlet_soft, col sep=comma]{example1_cond.csv};
      \label{pgfplots:dirichlet_soft}
      \addplot [very thick,cyan!70!black,mark=x,every mark/.append style={solid}]table[x=h, y=lacour_soft, col sep=comma]{example1_cond.csv};
      \label{pgfplots:lacour_soft}
      \addplot [very thick,orange!70!black,mark=triangle*,every mark/.append style={solid}]table[x=h, y=feti_soft, col sep=comma]{example1_cond.csv};
      \label{pgfplots:feti_soft}

      \nextgroupplot[align=left, title={{\scriptsize Hard inclusion}},ylabel={},yticklabels={}]
      \addplot [very thick,red!70!black,mark=diamond*,every mark/.append style={solid}]table[x=h, y=s_hard, col sep=comma]{example1_cond.csv};
\addplot [very thick,green!70!black,mark=square*,every mark/.append style={solid}]table[x=h, y=simple_hard, col sep=comma]{example1_cond.csv};
\addplot [very thick,blue!70!black,mark=*,every mark/.append style={solid}]table[x=h, y=dirichlet_hard, col sep=comma]{example1_cond.csv};
\addplot [very thick,cyan!70!black,mark=x,every mark/.append style={solid}]table[x=h, y=lacour_hard, col sep=comma]{example1_cond.csv};
\addplot [very thick,orange!70!black,mark=triangle*,every mark/.append style={solid}]table[x=h, y=feti_hard, col sep=comma]{example1_cond.csv};

      \nextgroupplot[align=left, title={{\scriptsize Uniform inclusion}},ylabel={},yticklabels={}]
      \addplot [very thick,red!70!black,mark=diamond*,every mark/.append style={solid}]table[x=h, y=s_equal, col sep=comma]{example1_cond.csv};
\addplot [very thick,green!70!black,mark=square*,every mark/.append style={solid}]table[x=h, y=simple_equal, col sep=comma]{example1_cond.csv};
\addplot [very thick,blue!70!black,mark=*,every mark/.append style={solid}]table[x=h, y=dirichlet_equal, col sep=comma]{example1_cond.csv};
\addplot [very thick,cyan!70!black,mark=x,every mark/.append style={solid}]table[x=h, y=lacour_equal, col sep=comma]{example1_cond.csv};
\addplot [very thick,orange!70!black,mark=triangle*,every mark/.append style={solid}]table[x=h, y=feti_equal, col sep=comma]{example1_cond.csv};
\end{groupplot}
    \matrix [ draw, matrix of nodes, anchor = west, node font=\scriptsize,
      column 1/.style={nodes={align=center,text width=0.5 cm}},
      column 2/.style={nodes={align=center,text width=0.6 cm}},
    ] at ($(group c3r1)+(2.2,0)$)
    {
      $\matI$                  & \ref{pgfplots:s_soft} \\
      $\matP_{S}$ & \ref{pgfplots:simple_soft} \\
      $\matP_{D}$ & \ref{pgfplots:dirichlet_soft} \\
      $\matP_{L}$ & \ref{pgfplots:lacour_soft} \\
      $\matP_{F}$ & \ref{pgfplots:feti_soft} \\
    };
\end{tikzpicture}
  \caption{Condition number of the preconditioned Schur complement system with respect to the mesh size.}
  \label{fig:convergence_1_cond}
\end{figure}

Now, we compare the performance of the preconditioned-CG (PCG) method for solving the problem \eqref{eq:saddle_AL}.
The termination criterion for the dual problem is chosen as $\|\vecy^{k+1}-\vecy^{k}\|_{\matS} < 10^{-12}$, and for the primal problem it is chosen as $\|\vecx^{k+1}-\vecx^{k}\|_{\matA} < 10^{-14}$.
In this comparison, we use the multigrid method as a preconditioner for the primal system.
For the experiments, we begin with solving the problem on level $L2$ with two-grid method and then increase the refinement level and also add another level in the multigrid hierarchy.
From Table~\ref{tab:example_iters}, we can observe that the FETI preconditioners are quite robust compared to the SIMPLE preconditioner.
As discussed earlier, SIMPLE preconditioner is only effective for uniform inclusion test case, while FETI preconditioners are stable irrespective of the mesh size and the material parameters.
The number of iterations stays constant for solving both the primal and dual systems.
In addition, the average number of iterations for solving the primal problem remain stable for a given test case.
Thus, the FETI preconditioner $\matP_F$  and the multigrid preconditioner are ideal for solving problem with discontinuities.
\begin{table*}[t]
  \begin{subtable}[]{1\linewidth}\centering
    \begin{tabular}{ c || c c c| c c c| c c c | c c c }
levels & $\matI$ & MG  & avg  & $\matP_S$ & MG  & avg  & $\matP_L$ & MG  & avg  & $\matP_F$ & MG  & avg  \\\hline\hline
      $L2$   & 23      & 216 & 9.39 & 31        & 285 & 9.19 & 16        & 153 & 9.36 & 11        & 108 & 9.82 \\\hline
      $L3$   & 24      & 225 & 9.38 & 33        & 306 & 9.27 & 17        & 162 & 9.53 & 11        & 108 & 9.82 \\\hline
      $L4$   & 28      & 263 & 9.39 & 34        & 322 & 9.47 & 16        & 155 & 9.69 & 11        & 109 & 9.91 \\\hline
      $L5$   & 28      & 261 & 9.32 & 34        & 315 & 9.26 & 17        & 162 & 9.53 & 11        & 108 & 9.82
    \end{tabular}
    \caption{Soft inclusion}
  \end{subtable}
  \begin{subtable}[]{1\linewidth}\centering
    \begin{tabular}{ c || c c c| c c c | c c c | c c c }
levels & $\matI$ & MG  & avg  & $\matP_S$ & MG  & avg  & $\matP_L$ & MG  & avg  & $\matP_F$ & MG  & avg  \\\hline\hline
      $L2$   & 24      & 201 & 8.38 & 33        & 273 & 8.27 & 16        & 137 & 8.56 & 11        & 97  & 8.82 \\\hline
      $L3$   & 24      & 225 & 9.38 & 32        & 297 & 9.28 & 17        & 162 & 9.53 & 12        & 117 & 9.75 \\\hline
      $L4$   & 29      & 270 & 9.31 & 34        & 315 & 9.26 & 17        & 162 & 9.53 & 11        & 108 & 9.82 \\\hline
      $L5$   & 28      & 262 & 9.36 & 35        & 325 & 9.29 & 17        & 163 & 9.59 & 12        & 118 & 9.83
    \end{tabular}
    \caption{Hard inclusion}
  \end{subtable}
  \begin{subtable}[]{1\linewidth}\centering
    \begin{tabular}{ c || c c c| c c c | c c c| c c c  }
levels & $\matI$ & MG  & avg  & $\matP_S$ & MG  & avg  & $\matP_L$ & MG  & avg  & $\matP_F$ & MG & avg  \\\hline\hline
      $L2$   & 35      & 240 & 6.86 & 23        & 153 & 6.65 & 15        & 99  & 6.60 & 11        & 78 & 7.09 \\\hline
      $L3$   & 36      & 229 & 6.36 & 24        & 163 & 6.79 & 15        & 101 & 6.73 & 12        & 87 & 7.25 \\\hline
      $L4$   & 43      & 265 & 6.16 & 25        & 163 & 6.52 & 14        & 91  & 6.50 & 12        & 81 & 6.75 \\\hline
      $L5$   & 43      & 265 & 6.16 & 25        & 161 & 6.44 & 15        & 97  & 6.47 & 11        & 74 & 6.73
    \end{tabular}
    \caption{Uniform inclusion}
  \end{subtable}
\caption{The number of PCG iterations required to reach a predefined tolerance, using various dual preconditioners and MG preconditioner, and average primal iterations per dual iteration.}
  \label{tab:example_iters}
\end{table*}

\subsection{Multiple inclusions}
\begin{figure*}[!htb]
  \centering
  \subcaptionbox{5 inclusions}[.24\textwidth]{\includegraphics{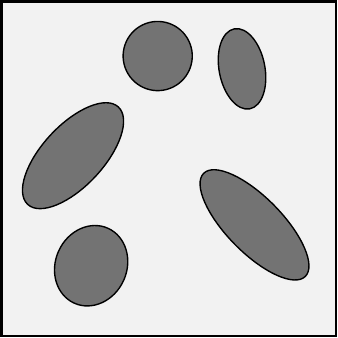}
  }
  \hfill \subcaptionbox{15 inclusions}[.24\textwidth]{\includegraphics{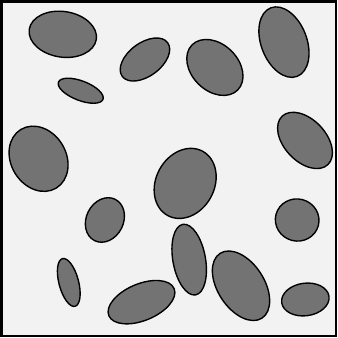}
  }
  \hfill
  \subcaptionbox{25 inclusions}[.24\textwidth]{\includegraphics{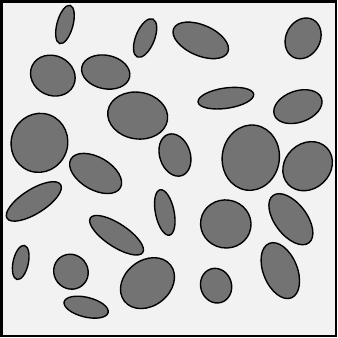}
  }
  \hfill
  \subcaptionbox{35 inclusions}[.24\textwidth]{
    \includegraphics{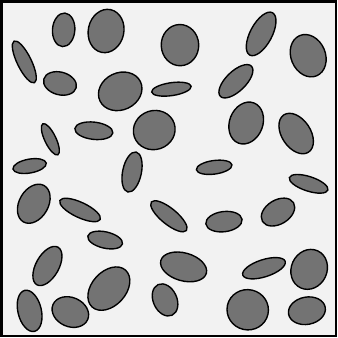}
  }
  \caption[2D mesh]{Sketch of domains with multiple inclusions.}
  \label{fig:multiple_inclusions}
\end{figure*}

\begin{figure*}[!htb]
  \centering
  \subcaptionbox{5 inclusions}[.24\textwidth]{\includegraphics[width=0.24\textwidth]{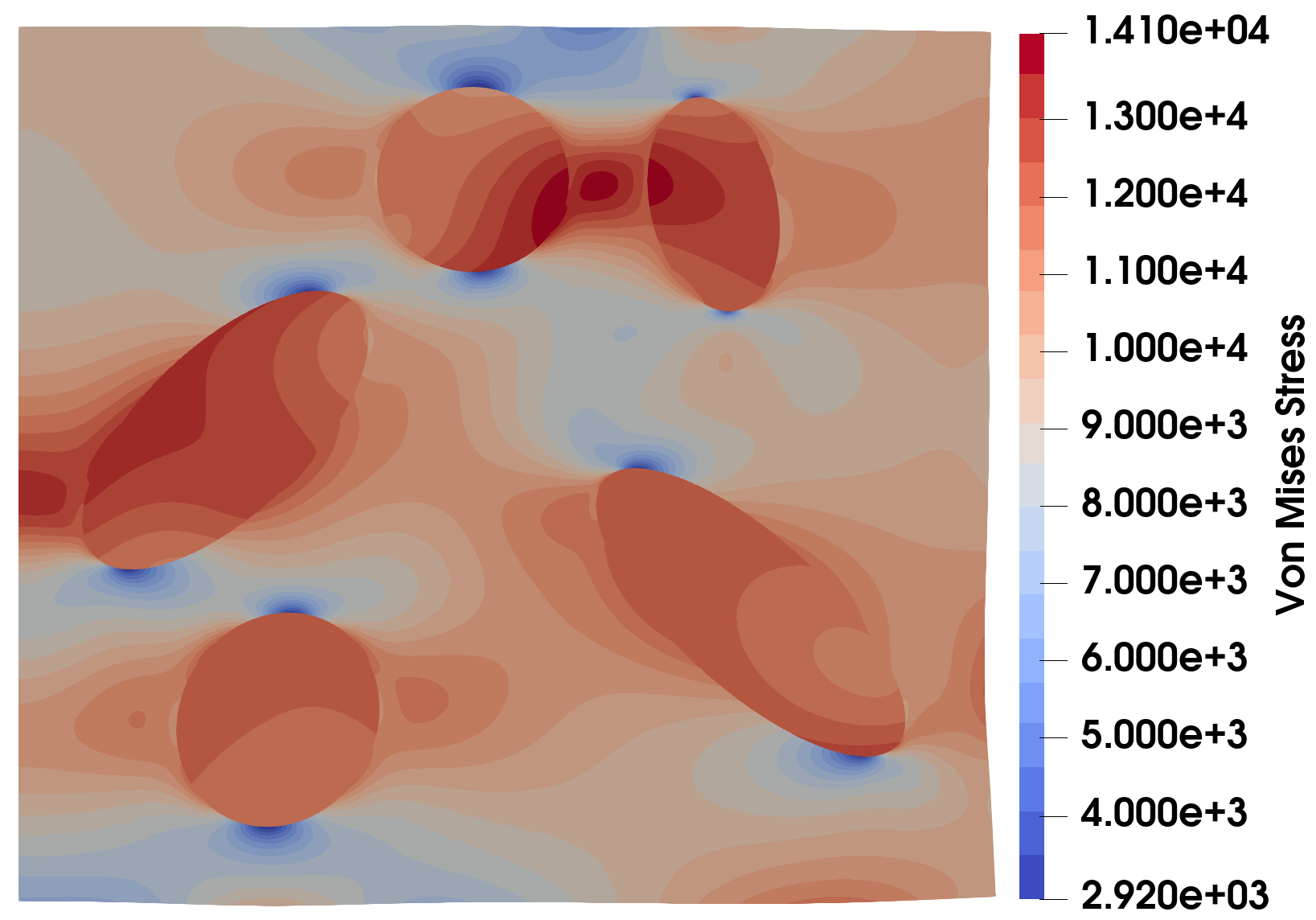}
  }
  \hfill \subcaptionbox{15 inclusions}[.24\textwidth]{\includegraphics[width=.24\textwidth]{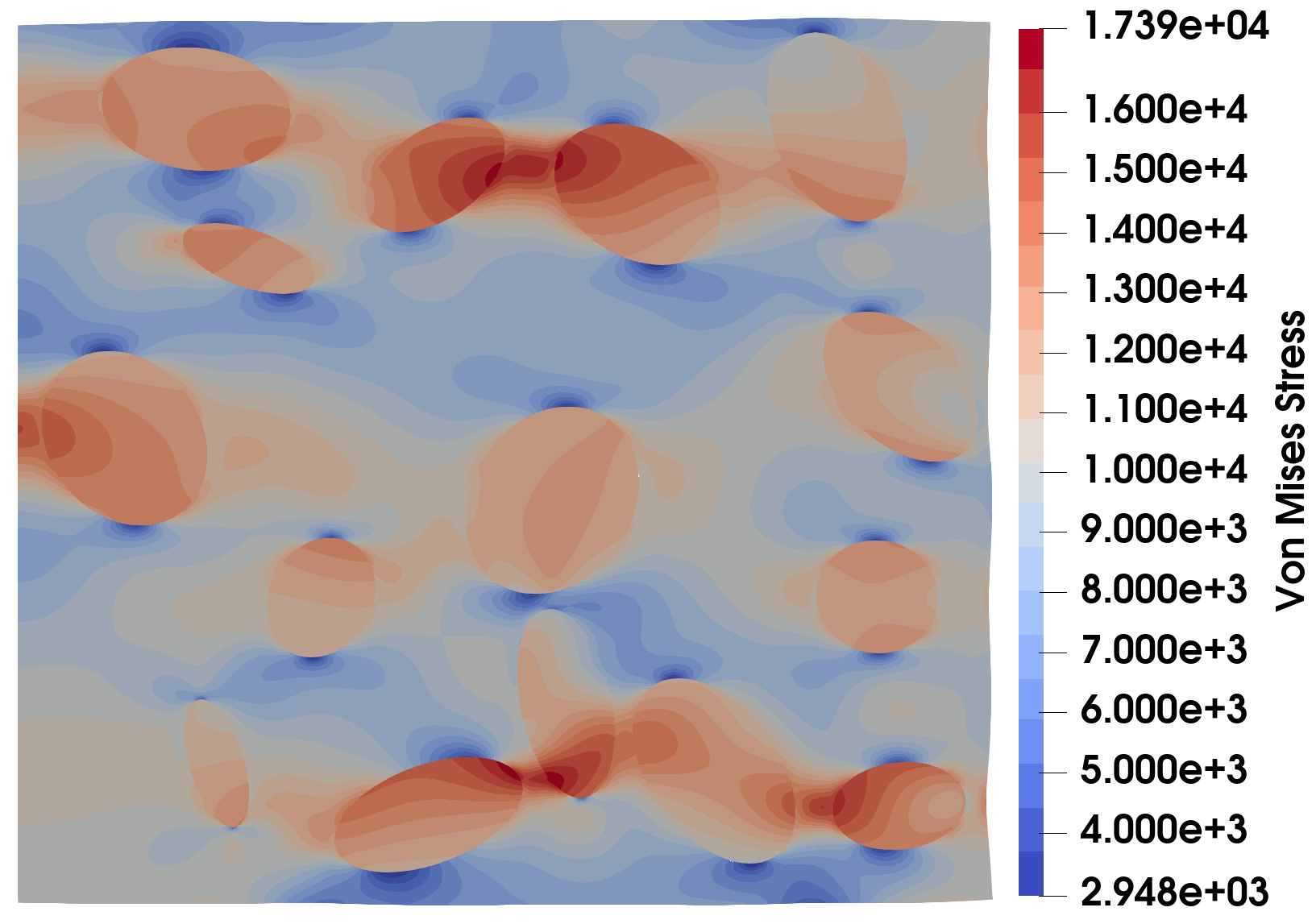}
  }
  \hfill
  \subcaptionbox{25 inclusions}[.24\textwidth]{\includegraphics[width=.24\textwidth]{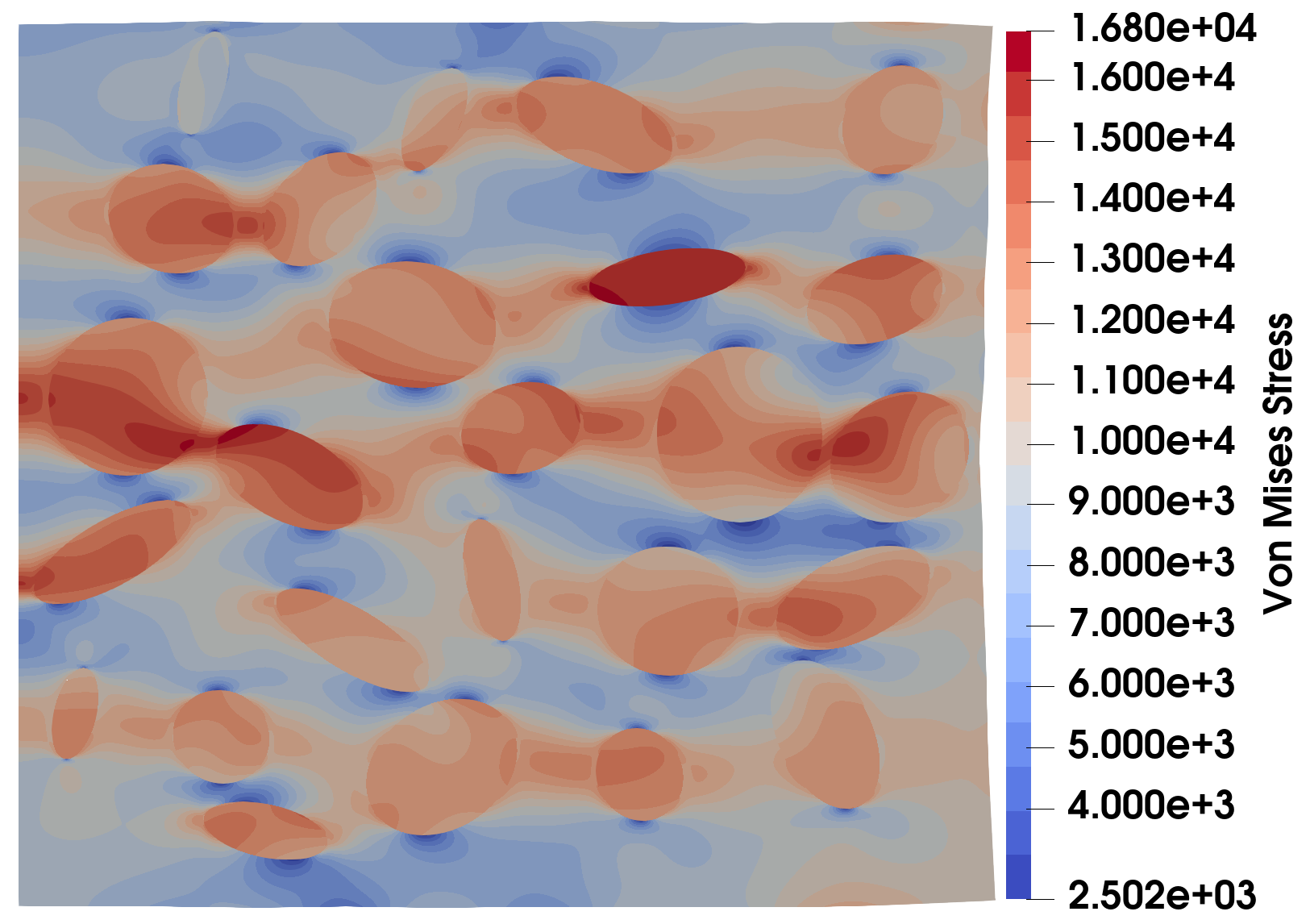}
  }
  \hfill
  \subcaptionbox{35 inclusions}[.24\textwidth]{
    \includegraphics[width=.24\textwidth]{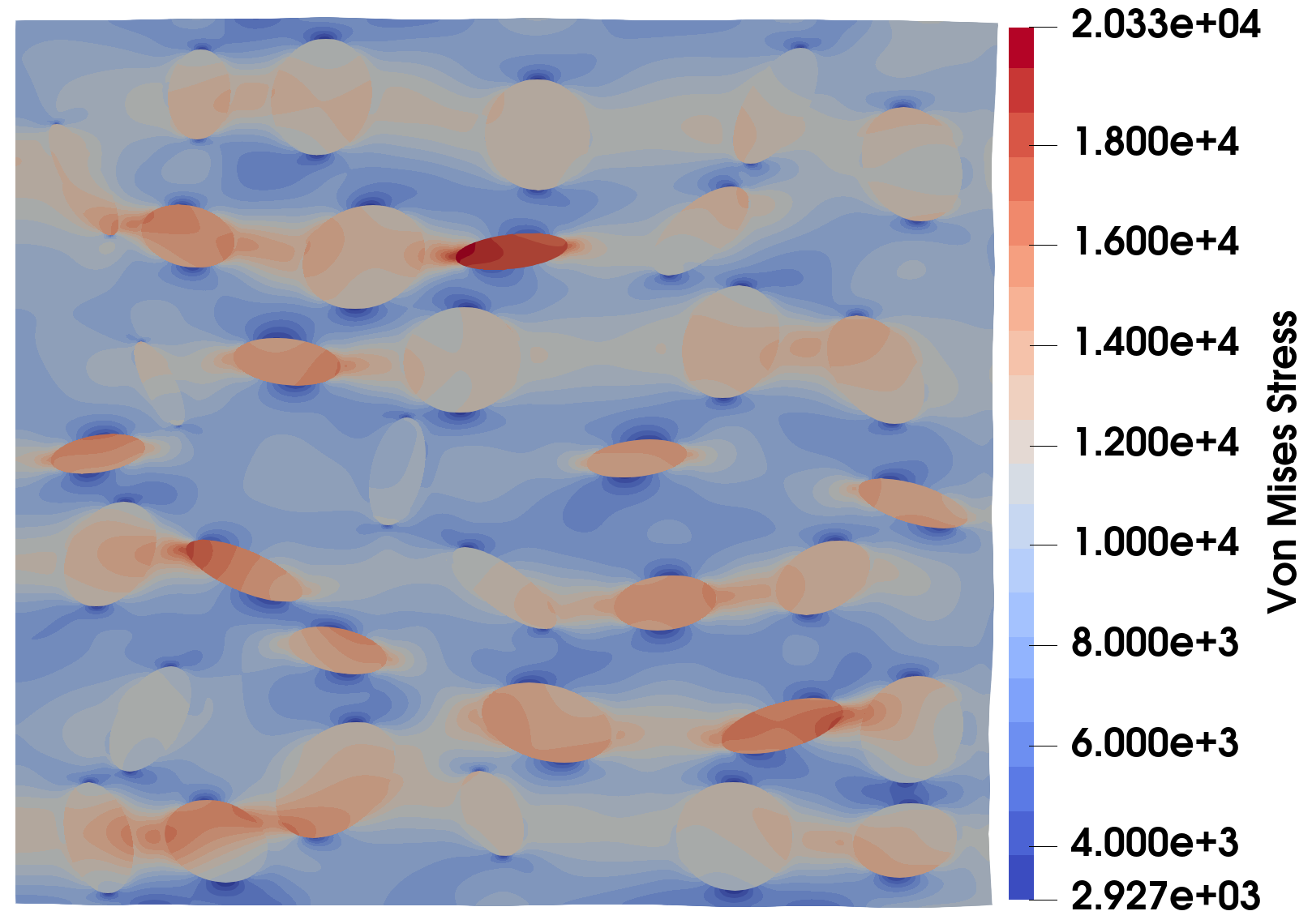}
  }
  \caption[2D mesh]{Induced Von Mises stress in a domain with multiple inclusions.}
  \label{fig:mulitple_result}
\end{figure*}

The second set of numerical experiments is carried out on a domain with multiple inclusions, which is subjected to the same boundary conditions as for the first experiments.
We randomly generate up to $35$ elliptical inclusions, which have uniformly distributed semi-major, semi-minor axes, orientation and the center~\cite{hiriyur_developments_2012}, as we can see in Figure~\ref{fig:multiple_inclusions}.
We define Young's modulus for inclusions $E_I=4\times 10^5\,\text{MPa}$ and the matrix $E_M = 10^5\,\text{MPa}$, while the Poisson's ration is $\nu=0.3$.
Figure~\ref{fig:mulitple_result} shows the resultant von Mises stress in the domain, since Young's modulus is higher for the inclusions, we can also see the induced stress is higher in the inclusions than in the matrix.

\begin{table*}[!htb]
  \centering
  \begin{tabular}{ c || c c c| c c c| c c c | c c c }
    \multirow{2}{*}{levels}            &
    \multicolumn{3}{c|}{5 inclusions}  &
    \multicolumn{3}{c|}{15 inclusions} &
    \multicolumn{3}{c|}{25 inclusions} &
    \multicolumn{3}{c} {35 inclusions}                                                                        \\
    \cline{2-13}
          & $\matP_F$ & MG  & avg  & $\matP_F$ & MG  & avg  & $\matP_F$ & MG  & avg  & $\matP_F$ & MG  & avg  \\\hline\hline
    $L2$  & 12        & 105 & 8.75 & 12        & 110 & 9.17 & 12        & 109 & 9.08 & 12        & 118 & 9.83 \\\hline
    $L3$  & 12        & 105 & 8.75 & 12        & 106 & 8.83 & 12        & 105 & 8.75 & 12        & 117 & 9.75 \\\hline
    $L4$  & 12        & 105 & 8.75 & 12        & 106 & 8.83 & 12        & 106 & 8.83 & 12        & 114 & 9.50 \\\hline
    $L5$  & 12        & 105 & 8.75 & 12        & 106 & 8.83 & 12        & 105 & 8.75 & 12        & 106 & 8.83
  \end{tabular}
  \caption{The number of PCG iterations required to reach a predefined tolerance, using a FETI dual preconditioners and MG preconditioner, and average primal iterations per dual iteration.}
  \label{tab:multiple_iters}
\end{table*}

In the previous section, we concluded that the FETI preconditioner $\matP_F$ outperformed the other considered preconditioners.
Hence, we employ the multigrid and the FETI preconditioner $\matP_F$ for solving the primal and dual systems, respectively.
Table~\ref{tab:multiple_iters} shows the number of required iterations for the convergence of both primal and dual systems.
We can see that the number of iterations for primal and dual systems are stable with respect to the mesh size and number of inclusions in the domain.
 
\section{Conclusions}
In this paper, we discussed an unfitted FE method for modeling the inclusions.
We demonstrated that the method has optimal convergence properties for different test cases.
We discussed the SIMPLE and FETI preconditioners for solving the dual systems, where the FETI preconditioner $\matP_F$ outperformed the other preconditioners.
For solving the primal system, we discussed a multigrid method that demonstrates the level-independence property.
We show that the proposed preconditioners are stable and robust, as the number of iterations for both primal and dual systems does not vary with respect to the number of inclusions in the domain.

\bibliographystyle{spmpsci}      \def\url#1{}

\end{document}